\documentstyle[amssymb,amstex,amsopn,graphicx,graphics,12pt,a4wide]{article}

\newtheorem{theorem}{Theorem}

\newtheorem{proposition}{Proposition}

\newtheorem{remark}{Remark}
\newtheorem{example}[theorem]{Example}
\def\EE{{\mathbb{E}}}
\def\PP{{\mathbb{P}}}
\def\HH{{\mathbb{H}}}
\def\cB{{\mathcal{B}}}
\def\cV{{\mathcal{V}}}
\def\cZ{{\mathcal{Z}}}
\def\DH{\text{dim}_H}
\def\CD{\text{cardim\ }}
\def\diam{\text{diam}}
\def\len{\text{length}}
\def\DM{\overline{\text{dim}}_M}
\DeclareMathOperator{\esssup}{ess\,sup}
\DeclareMathOperator{\essinf}{ess\,inf} 

\begin{document}

\title{\textbf{Hausdorff dimension of visibility sets for well-behaved continuum percolation in the hyperbolic plane}}
\author{Christoph Th\"ale\footnote{\textit{e-mail: christoph.thaele[at]uni-osnabrueck.de}}\\
Department of Mathematics\\
University of Osnabr\"uck, Germany}
\date{}
\maketitle

\begin{abstract} Let $\cZ$ be a so-called well-behaved percolation, i.e. a certain random closed set in the hyperbolic plane, whose law is invariant under all isometries; for example the covered region in a Poisson Boolean model. The Hausdorff-dimension of the set of directions is determined in terms of the $\alpha$-value of $\cZ$ in which visibility from a fixed point to the ideal boundary of the hyperbolic plane is possible within $\cZ$. Moreover, the Hausdorff-dimension of the set of (hyperbolic) lines through a fixed point contained in $\cZ$ is calculated. Thereby several conjectures raised by Benjamini, Jonasson, Schramm and Tykesson are confirmed.
\end{abstract}
\begin{flushleft}\footnotesize
\end{flushleft}

\section{Introduction and main result}

In this note we are interested in some well-behaved percolation models in the hyperbolic plane. This topic has been considered by several authors and became an active field of current research, see \cite{BJST,BenjaminiScramm,CT,Lalley,Tyk} to name just a few. Background material on hyperbolic geometry may be found in \cite{BenePetrino,Ramseyetal} and some aspects of percolation theory in the Euclidean spaces is presented in \cite{MR}. Our focus here will be on the set of hyperbolic lines (bi-infinite geodesic rays) and half-lines (infinite geodesic rays) contained in the unbounded connected components of some continuum percolation models in the hyperbolic plane. Of course, similar problems can also be treated in higher dimensional hyperbolic spaces $\HH^d$ or the $d$-dimensional Euclidean space ${\Bbb R}^d$. However, it has been shown (see \cite{BJST,VisibilityPierre}) that for example $2$-dimensional planes that are contained in $\cZ$ do not exist for well behaved-percolation in ${\Bbb H}^d$ for any $d\geq 3$. Moreover, visibility to infinite in ${\Bbb R}^d$ is impossible even for $2\geq 2$. For this reason we restrict our attention to the hyperbolic plane $\HH^2$ equipped with the usual hyperbolic standard metric $\varrho_{\HH^2}$.

To formulate our main results, which confirm several conjectures raised by Benjamini, Jonasson, Schramm and Tykesson in \cite{BJST}, let $B(1)\subset\HH$ be a closed disc of radius $1$. A random closed set $\cZ$ in $\HH^2$ is called a \textbf{well-behaved percolation} if the following assumptions are satisfied (see \cite{BJST}):
\begin{itemize}
 \item[(i)] The law of $\cZ$ is invariant under all isometries of $\HH^2$.
 \item[(ii)] For any two bounded increasing measurable functions $g$ and $h$ of $\cZ$, the FKG-type inequality $$\EE[g(\cZ)h(\cZ)]\geq\EE[g(\cZ)]\EE[h(\cZ)]$$ is satisfied.
 \item[(iii)] There is some $R_0<\infty$ such that $\cZ$ satisfies independence at distance $R_0$. This is for every subsets $A,B\subset\HH^2$ with $\inf\{\varrho_{\HH^2}(a,b):a\in A,b\in B\}\geq R_0$ the events $\cZ\cap A$ and $\cZ\cap B$ are independent.
 \item[(iv)] The expected number of connected components of $B(1)\setminus\cZ$ is finite.
 \item[(v)] We have $\EE[\len(B(1)\cap\partial\cZ)]<\infty$.
 \item[(vi)] We have $\PP(B(1)\subset\cZ)>0$.
\end{itemize}
Let $f(r)$ denote the probability that a fixed line segment of length $r>0$ is contained in $\cZ$ and fix some point $o\in\HH^2$. We recall from \cite{BJST}, Lemma 3.4 that there exists a unique $\alpha\geq 0$ called the \textbf{$\alpha$-value of $\cZ$} such that $f(r)=\Theta(e^{-\alpha r})$ for any $r\geq 0$ in the usual Landau notation. In terms of its $\alpha$-value, the Hausdorff-dimension of several random sets related to a well-behaved percolation $\cZ$ in $\HH^2$ can be determined: 
\begin{theorem}\label{thm1} Consider a well-behaved percolation $\cZ$ in $\HH^2$ and a fixed point $o\in\HH^2$. Let $\frak V$ denote the set of points $z$ in the ideal boundary $\partial\HH^2$ of the hyperbolic plane such that the ray $[o,z)$ is contained in $\cZ$. If $\alpha\geq 1$ then ${\frak V}=\emptyset$ with probability one. If $\alpha<1$ then $\PP({\frak V}\neq\emptyset)>0$ and $\DH {\frak V}=1-\alpha$ almost surely on ${\frak V}\neq\emptyset$. Moreover, the union of all these rays has Hausdorff-dimension $2-\alpha$ almost surely on ${\frak V}\neq\emptyset$.
\end{theorem}
\begin{theorem}\label{thm2}
 For a well-behaved percolation $\cZ$ in the hyperbolic plane $\HH^2$ and fixed $o\in\HH^2$ we have: If $\alpha\geq 1/2$ then there is no line through $o$ contained in $\cZ$ almost surely. If $\alpha<1/2$ then the union of all lines in $\cZ$ through $o$ has Hausdorff-dimension $2-2\alpha$ with probability one conditioned on the event that there are such lines.
\end{theorem}
The following random sets are \textbf{examples} to which our theory applies (see \cite{BJST,CT}):
\begin{example} Let $\eta_\lambda$ be an isometry-invariant Poisson point process of intensity $\lambda\in(0,\infty)$ in $\HH^2$, $R>0$ and define $$\cB:=\bigcup_{x\in\eta_\lambda}B(x,R)\ \ \ \text{and}\ \ \ \cV:=\overline{\HH^2\setminus\cB}.$$ Then, $\cB$ and $\cV$, the occupied and the vacant phase of the Boolean model with respect to $\eta_\lambda$ and $R$, are well-behaved percolation sets. The $\alpha$-value for $\cV$ is given by $\alpha=2\lambda\sinh R$ and the $\alpha$-value for $\cB$ is the unique solution of $$\int_0^{2R}e^{\alpha t}H_{\lambda,R}(t)dt=1,$$ where $$H_{\lambda,R}(t)=-\exp\left(-4\lambda\int_0^{t\over 2}\sinh\left(\cosh^{-1}\left({\cosh R\over\cosh s}\right)\right)ds\right).$$
\end{example}
\begin{example}
Let $\eta_\lambda$ be as above but consider the radius $R$ of the balls in the Boolean model as random and assume that the exponential moment $\EE[e^R]$ is finite. In this case the $\alpha$-value of the vacant phase $\cV$ equals $\alpha=2\lambda\EE[\sinh R]$.
\end{example}
\begin{example}
Let $\eta_\lambda$ be as in our previous examples and let $K$ be a random closed convex set with a.s. finite diameter containing the origin, whose law is invariant under isometries of $\HH^2$. We can think of $\HH^2$ as the unit disc embedded in the complex plane and put $\varphi_x(z)=(z-x)/(1-\overline{x}z)$, where $\overline{\cdot}$ stands for complex conjugation. Let us define $$\cB_K=\bigcup_{x\in\eta}\varphi_x^{-1}(K_x)\ \ \ \text{and}\ \ \ \cV_K=\overline{\HH^2\setminus\cB},$$ where $\{K_x:x\in\eta\}$ is an i.i.d. family of random sets indexed by the points of $\eta_\lambda$ having the same distribution as $K$ and are independent of $\eta_\lambda$. Then $$f(r)=\exp(-\lambda\EE[\text{Area}(x\in\HH^2:\varphi_x^{-1}(K)\cap L_{o,r}\neq\emptyset)])$$ for the well-behaved percolation $\cV_K$, where $L_{o,r}$ is a line segment of length $r>0$ starting at $o$. Thus, $\alpha=\EE[\text{Area}(x\in\HH^2:\varphi_x^{-1}(K)\cap L_{o,r}\neq\emptyset)]$.
\end{example} 
The rest of this note is organized as follows: In Section \ref{secaux} we recall some facts from fractal geometry and geometric measure theory and prove an auxiliary result on Hausdorff-dimensions of random sets. In the final section we present the proofs of our main results.

\section{An auxiliary result on Hausdorff dimensions of random sets}\label{secaux}

Let $(E,\varrho)$ be a metric space, which is second countable, locally compact and has the Hausdorff property (a so-called lcscH space). Let ${\cal B}$ be the Borel $\sigma$-field on $E$ generated by $\varrho$, ${\cal F}$ be the family of closed subsets of $E$ and let ${\cal M}$ be the family of Radon measures on $E$ (recall that a \textbf{Radon measure} is a locally finite and inner regular measure on $\cal B$). We equip $\cal F$ with the $\sigma$-field $\frak F$ generated by the usual Fell-topology \cite[Appendix B]{Mol} on ${\cal F}$ and ${\cal M}$ with the $\sigma$-field $\frak M$ generated by the evaluation mappings $\varphi\mapsto\varphi(B)$, $B\in{{\cal B}}$, $\varphi\in{{\cal M}}$ (cf. Chapter 1.1 \cite{Kallenberg}). For $D\geq 0$ and $B\subset E$ the \textbf{$D$-dimensional Hausdorff-measure} ${{\cal H}}^D(B)$ is defined by $${\cal H}^D(B):=\lim_{\delta \downarrow 0}{\cal H}_\delta^D(B),$$ where $${\cal H}_\delta^D(B):=\inf_{{\cal V}}\left\lbrace\sum_{F\in{{\cal V}}}\alpha(d)(\diam(F))^D:B\subseteq\bigcup_{F\in{{\cal V}}}F,\ \diam(F)<\delta,\ F\in{{\cal F}}\right\rbrace$$ and where the infimum is taken over all countable subfamilies ${\cal V}$ of ${\cal F}$. Moreover, we put $\alpha(D)=\Gamma(1/2)^D/(2^D\Gamma(1+D/2))$. For $B\subset E$ the \textbf{Hausdorff-dimension} $\DH B$ of $B$ is defined by $$\DH B=\inf\left\lbrace D\geq 0:{{\cal H}}^D(B)=0\right\rbrace=\sup\left\lbrace D\geq 0:{{\cal H}}^D(B)=+\infty\right\rbrace.$$ For $D\geq 0$, the \textbf{${{\cal H}}^D$-derivative} of $\varphi\in{{\cal M}}$ at $x\in E$ is given via
\begin{eqnarray}
\nonumber {{\cal D}}(\varphi,D,x) &:=& \limsup_{F\rightarrow x}{\varphi(F)\over\alpha(D)(\diam(F))^D}\\
\nonumber &=& {1\over\alpha(D)}\limsup_{\delta\downarrow 0}\left\lbrace {\varphi(F)\over(\diam(F))^D}:x\in F,\ F\in{{\cal F}}, \diam(F)\leq\delta\right\rbrace. 
\end{eqnarray}
Denote $E^{(\infty)}=E^{(\infty)}(\varphi,D)=\{x\in E:{{\cal D}}(\varphi,D,x)=+\infty\}$.\\ A \textbf{random measure} $\eta$ on $E$ is a $[{{\cal M}},{\frak M}]$-valued random variable defined on some abstract probability space, cf. \cite{Kallenberg}. Its \textbf{second-moment measure} $\Lambda=\Lambda_\eta$ on $E\times E$ is defined by the relation
\begin{equation}\label{smm}
\Lambda(B\times B'):=\EE[\eta(B)\eta(B')],\ \ \ \ \ B,B'\in{{\cal B}}.
\end{equation}
We are now in the position to rephrase a Frostman-type result, which was proved in \cite{ZCutout,Zcharater} for the special case $E={\Bbb R}^n$. For completeness and to keep the argument below self-contained we include a streamlined proof in our more general setting. Later on the result will be applied to subsets of the hyperbolic plane.
\begin{proposition}\label{prophdim}
Let $\eta$ be a random measure on $E$, $D\geq 0$ and $r>0$. Suppose there exist a sequence $E_n\uparrow E$ with $E_n\in{{\cal B}}$ satisfying $$\int\varrho(x,y)^{-D}{\bf 1}[x\in E_n]{\bf 1}[\varrho(x,y)<r]\Lambda(d(x,y))<\infty$$ for any $n\in{\Bbb N}$. Then for $B\in{{\cal B}}$, almost surely on $\eta(B)>0$ we have $\DH B\geq D$.
\end{proposition}
\paragraph{Proof.} The proof is divided into four steps. 
\subparagraph{\textit{Step 1:}} \textit{If $\varphi\in{{\cal M}}$, then the restriction $\varphi\llcorner(E\setminus E^{(\infty)})$ of $\varphi$ onto $E\setminus E^{(infty)}$ is absolutly continuous with respect to ${{\cal H}}^D$.}\\ Define $\psi:=\varphi\llcorner(E\setminus E^{(\infty)})$, put for $a>0$, $E^{(a)}:=\{x\in E:{{\cal D}}(\varphi,D,x)\in[0,a)\}$ and let $B\subset E$ with ${{\cal H}}^D(B)=0$. Then 2.10.17 (3) in \cite{Federer} implies
\begin{eqnarray}
\nonumber \psi(B) &\leq& \psi(B\cap E^{(\infty)})+\psi(B\setminus E^{(\infty)})\leq 0+\lim_{a\rightarrow\infty}\varphi(B\cap E^{(a)})\\
\nonumber &\leq& \lim_{a\rightarrow\infty}a{{\cal H}}^D(B\cap E^{(a)})\leq\lim_{a\rightarrow\infty}a{{\cal H}}^D(B)=0.
\end{eqnarray}
\subparagraph{\textit{Step 2:}} \textit{If $B\subset E$, $\varphi\in{{\cal M}}$ and $\varphi(B\setminus E^{(\infty)})>0$, then $\DH B\geq D$.}\\ Indeed, using Step 1 we see that $\varphi(B\setminus E^{(\infty)})>0$ implies $0<{{\cal H}}^D(B\setminus E^{(\infty)})\leq{{\cal H}}^D(B)$. Thus, $\DH B\geq D$ by the definition of Hausdorff-dimension.
\subparagraph{\textit{Step 3:}} \textit{If $D\geq 0$, $B\subset E$, $\varphi\in{{\cal M}}$ with $\varphi(B)>0$ and for $\varphi$-almost all $x\in B$ there exists $r=r(x)>0$ with $\int_{B(x,r)}\varrho(x,z)^{-D}\varphi(dz)<\infty$, where $B(x,r)$ is the ball of radius $r$ around $x$, then $\DH B\geq D$.}\\ To see it, note that for $\varphi$-almost all $x\in B$ and any $F\in{{\cal F}}$ with $x\in F$ we have
\begin{eqnarray}
\nonumber \alpha(D){{\cal D}}(\varphi,D,x) &=& \limsup_{F\rightarrow x}{\varphi(F)\over(\diam(F))^D}\leq\limsup_{F\rightarrow x}\int_{F}\varrho(x,z)^{-D}\varphi(dz)\\
\nonumber &\leq& \int_{B(x,r)}\varrho(x,z)^{-D}\varphi(dz)<\infty.
\end{eqnarray}
Hence, $\varphi(B\cap E^{(\infty)})=0$ and $\varphi(B\setminus E^{(\infty)})=\varphi(B)>0$ and Step 2 implies $\DH B\geq D$.
\subparagraph{\textit{Step 4:}} Use Campbell's theorem to conclude $$\hspace{-5cm}{\Bbb E}\int_{E_n}\int_{B(x,r)}\varrho(x,y)^{-D}\eta(dy)\eta(dx)$$ $$\hspace{2cm}=\int\varrho(x,y)^{-D}{\bf 1}[x\in E_n]{\bf 1}[\varrho(x,y)<r]\Lambda(d(x,y))<\infty.$$ This implies that for $P_\eta$-almost all $\eta$ (here $P_\eta$ is the distribution of $\eta$) and $\eta$-almost all $x\in\bigcup E_n=E$ we have $$\int_{B(x,r)}\varrho(x,y)^{-D}\eta(dy)<\infty.$$ By Step 3 we see now that on $\eta(B)>0$ we have $\DH B\geq D$ almost surely.\hfill $\Box$ \\ \\ Let us further recall that a \textbf{random closed set} on the metric space $(E,\varrho)$ is a measurable mapping from some abstract probability space into the measurable space $[{\cal F},{\frak F}]$, see \cite{Mol}.

\section{Proofs}

\subsection{Preliminaries}

We recall that $f(r)$ denotes the probability that a fixed line segment of length $r>0$ is contained in $\cZ$. Moreover, for some fixed point $o\in\HH^2$, $A$ stands for a closed half-plane with $o$ on its boundary. We define $I:=A\cap\partial B(o,1)$. Furthermore, for $r>1$, $Y_r$ is the set of those $x\in I$ with the property that the line segment with endpoint $o$ through $x$ having length $r$ is contained in $\cZ$. Moreover, the random set $Y$ is defined by $$Y=\bigcap_{n\geq 1}Y_{nR_0},$$ where, recall, $R_0$ is the independence distance from the definition of $\cZ$.

In Lemma 3.6 of \cite{BJST} the following has been shown:
\begin{proposition}\label{prophelp} Let $\cZ$ be a well-behaved percolation in $\HH^2$. Then $f(r)\leq e^{-\alpha r}$, $$\EE[\len(Y_r)]=\len(I)f(r)$$ and $$\PP(x,y\in Y_r)\leq f(r)f(r+\log\varrho_{\HH^2}(x+y)+O(1)),\ \ \ x,y\in I.$$  
\end{proposition}
For the proof of Theorem \ref{thm1} we will estimate the Hausdorff-dimension $\DH S$ of a set $S$ from above by its upper Minkowski-dimension $\DM S$. In the case that the ambient space is the boundary $\partial C$ of a circle $C$ and $S\subset\partial C$, the upper Minkowski-dimension of $S$ relative to $\partial C$ can be defined by $$\DM S=1-\limsup_{\delta\rightarrow 0}{\log\len(S(\delta))\over\log\delta}.$$ Here, $S(\delta)$ stands for the $\delta$-parallel set of $S$ relative to $\partial C$, see \cite{Falconer}. Let us further recall the following well known inequality between Hausdorff- and upper Minkowski-dimension: $$\DH S\leq \DM S.$$

\subsection{Proof of Theorem \ref{thm1}} 

\paragraph{\textit{The case $\alpha\geq 1$:}} It has been shown in Lemma 3.5 of \cite{BJST} that for $\alpha\geq 1$ we have $\frak V=\emptyset$ with probability one. We can henceforth restrict our attention to the case $\alpha<1$, where the event ${\frak V}\neq\emptyset$ has positive probability.
\paragraph{\textit{An upper bound for the mean:}} We start by observing that $$\EE[\DM Y]\leq 1-\limsup_{r\rightarrow\infty}{\log\EE[\len(Y_{r})]\over -r}.$$ This is because parallel sets are taken in $I\subset B(o,1)$ and because $Y_r$ is an almost surely decreasing familiy of subsets of $I$. We now use Proposition \ref{prophelp}, which says that $\EE[\len(Y_r)]=\text{length}(I)f(r)$, thus, with $f(r)\leq e^{-\alpha r}$ and the inequality between Hausdorff- and Minkowski-dimension, we deduce 
\begin{eqnarray}
\nonumber \EE[\DH Y] &\leq& \EE[\DM Y] \leq 1-\limsup_{r\rightarrow\infty}{\log[\text{length}(I)f(r)]\over -r}\\
\nonumber &\leq& 1-\limsup_{r\rightarrow\infty}{\log[\text{length}(I)e^{-\alpha r}]\over -r}=1-\alpha.
\end{eqnarray}
\paragraph{\textit{A lower bound with positive probability:}} Let ${{\cal M}}_I$ be the space of Radon measure on $I$ equipped with the weak topology and define for $n\geq 1$ the random measure $\nu_n$ by $$d\nu_n:=e^{\alpha R_0n}{\bf 1}[\cdot\in Y_{R_0n}]dx,$$ where $dx$ stands for the element of the Lebesgue measure $I$ and where $R_0$ is the independence distance from the definition of well-behaved percolation. Obviously, $\nu_n\in{{\cal M}}_I$ and $||\nu_n||<\infty$ with probability one. Indeed, we have from Proposition \ref{prophelp}, $\EE[||\nu_n||]\leq\text{length}(I)<\infty$, which implies $||\nu_n||<\infty$ almost surely. Moreover, for any Borel set $B\subset I$ we have by Markov's inequality and Proposition \ref{prophelp}
\begin{eqnarray}
\nonumber \lim_{t\rightarrow\infty}\limsup_{n\rightarrow\infty}\PP(\nu_n(B)>t) &= & \lim_{t\rightarrow\infty}\limsup_{n\rightarrow\infty}\PP(e^{\alpha R_0n}\len(B\cap Y_{R_0n})>t)\\
\nonumber &\leq & \lim_{t\rightarrow\infty}\limsup_{n\rightarrow\infty}{e^{\alpha R_0n}\EE[\len(B\cap Y_r)]\over t}\\
\nonumber &\leq & \lim_{t\rightarrow\infty}\limsup_{n\rightarrow\infty}{e^{\alpha R_0n}\EE[\len(Y_{R_0n})]\over t}=0.
\end{eqnarray}
A similar argument also shows $$\inf_B\limsup_{n\rightarrow\infty}\PP(\nu_n(B^C)>\varepsilon)=0,\ \ \ \varepsilon>0,$$ where the infimum is taken over all Borel sets $B\subset I$. We can now apply Lemma 4.5 and 4.11 in \cite{Kallenberg} to conclude that the sequence $(\nu_n)$ is relatively compact with respect to the weak topology on ${\cal M}_I$. Thus, any sequence $(n)$ contains a subsequence $(n')$ such that $\nu_{n'}$ converges weakly to some limit measure, which is almost surely bounded. Moreover, the second-moment estimate $\EE[||\nu_n||^2]=O(1)(\EE[||\nu_n||])^2$ has been shown in \cite{BJST}. Thus, there exists $\varepsilon>0$ such that $\PP(||\nu_n||>\varepsilon)>0$ for all $n$. Hence, with positive probability we can extract a subsequence $\nu_{n_k}$, such that $||\nu_{n_k}||>\varepsilon$ for all $k$. Moreover, by a compactness argument we can pass to a further subsequence that converges weakly to some limit measure $\nu$ satisfying $||\nu||>0$. For this reason the measure $\nu$ can be regarded as a mass distribution on the intersection $Y=\bigcap_{n>1}Y_{R_0n}$, provided $Y$ is not empty.\\ We consider now the second-moment measure $\Lambda=\Lambda_\nu$ of $\nu$, which can be defined as in (\ref{smm}). We have by Fatou's lemma and Fubini's theorem
\begin{eqnarray}
\nonumber \Lambda(B\times B') &\leq& \lim_{n\rightarrow\infty}{\Bbb E}\int_B\int_{B'}e^{2\alpha R_0n}{\bf 1}[x,y\in Y_{R_0n}]dxdy\\
\nonumber &=& \lim_{n\rightarrow\infty}\int_B\int_{B'}e^{2\alpha R_0n}\PP(x,y\in Y_{R_0n})dxdy
\end{eqnarray}
for Borel sets $B,B'\subset I$. Furthermore, from the second-moment estimate in Proposition \ref{prophelp} it follows $$\PP(x,y\in Y_{R_0n})\leq f(R_0n)f(R_0n+\log\varrho_{\HH^2}(x,y)+O(1))\leq e^{-2\alpha R_0n}\varrho_{\HH^2}(x,y)^{-\alpha}O(1)$$ for any $n\geq 1$, whence $$\Lambda(B\times B')\leq O(1)\int_B\int_{B'}{dxdy\over\varrho_{\HH^2}(x,y)^\alpha}.$$ We now observe that $$\int_{Y}\int_{Y}{\Lambda(d(x,y))\over\varrho_{\HH^2}(x,y)^D}\leq O(1)\int_I\int_I{dxdy\over\varrho_{\HH^2}(x,y)^{D+\alpha}}$$ is finite whenever $D+\alpha<1$, or equivalently, if $D<1-\alpha$. Hence, together with Proposition \ref{prophdim} () we see that there is positive probability for the event $\DH Y\geq 1-\alpha$.
\paragraph{\textit{A lower bound with probability one:}} It remains to show that we have $\DH Y\geq 1-\alpha$ with probability one on $Y\neq\emptyset$. To this end denote by ${\cal F}_n$ the $\sigma$-field generated by $Y_{R_0n}$, observe that all these $\sigma$-fields are independent, because of the definition of $R_0$. Define further ${\cal A}_n$ as the $\sigma$-field generated by the family $\{{\cal F}_m:m\geq n\}$ and put ${\cal T}:=\bigcap_{n\geq 1}{\cal A}_n$. It is easily checked that $\{\DH Y\geq s\}\in{\cal T}$ for any $s\in[0,\infty)$. Thus the 0-1-law, Theorem 3.13 in \cite{Kallenberg2}, implies that the event $\{\DH Y\geq s\}$ has probability $0$ or $1$. On the other hand, we have shown that $\DH Y\geq 1-\alpha$ holds on $Y\neq\emptyset$ with positive probability, which allows us to conclude $\DH Y\geq 1-\alpha$ almost surely on $Y\neq\emptyset$.
\paragraph{\textit{The ideal boundary:}} So far we have proved that $$\EE[\DH Y]\leq 1-\alpha\ \ \ \text{and that}\ \ \ \PP(\DH Y\geq 1-\alpha|Y\neq\emptyset)=1,$$ which clearly implies $\DH Y=1-\alpha$ with probability one on $Y\neq\emptyset$. But this value is independent of the choice of the defining half-plane $A$, which implies by invariance of $\cZ$ that the random set $Y'\subset B(o,1)$ of those $x$ for which the hyperbolic half-line (ray) through $x$ starting at $o$ is fully contained in $\cZ$ has also Hausdorff-dimension $1-\alpha$ with probability one on $Y'\neq\emptyset$. However, this is obviously the same as the set $\frak V$ of points $z$ on the ideal boundary $\partial\HH^2$ for which $[o,z)\subset\cZ$, which proves $$\PP(\DH{\frak V}=1-\alpha|{\frak V}\neq\emptyset)=1.$$
\paragraph{\textit{The set of rays:}} We denote by ${\cal R}_o$ the set of hyperbolic rays $[o,z)$ with $z\in\partial\HH^2$ and the property that $[o,z)\subset\cZ$. Defining ${\cal R}_o':=Y'\times[0,1]$ with $Y'$ as in the previous paragraph, standard fractal geometry (see Corollary 7.4 in \cite{Falconer}) implies $$\DH{\cal R}_o'=\DH Y'+\DH[0,1]=2-\alpha$$ almost surely on $Y'\neq\emptyset$. It is readily verified that this implies $$\DH{\cal R}_o=\DH{\frak V}+1=2-\alpha\ \ \ \ \text{a.s. on}\ {\cal R}_o\neq\emptyset,$$ which finally completes the proof.\hfill $\Box$  

\begin{remark}
Let $\mu$ be a Radon measure defined on some nice metric space $E$ as in Section \ref{secaux} above and define the \textbf{lower and upper pointwise dimension of $\mu$ at $x\in E$} as $$\underline{\text{d}}\mu(x)=\liminf_{r\rightarrow 0}{\ln\mu(B(x,r))\over\ln r}\ \ \ \text{and}\ \ \ \overline{\text{d}}\mu(x)=\limsup_{r\rightarrow 0}{\ln\mu(B(x,r))\over\ln r},$$ respectively. Moreover, the \textbf{lower and upper Hausdorff-dimension of $\mu$} are given by $$\underline{\text{d}}_H\mu=\essinf_{x\in E}\underline{d}\mu(x)\ \ \ \text{and}\ \ \ \overline{\text{d}}_H\mu=\esssup_{x\in E}\overline{\text{d}}\mu(x),$$ respectively. The measure $\mu$ is saied to have \textbf{carrying dimension} $\beta$, this is $\CD\mu=\beta$, if $\underline{\text{d}}_H\mu=\overline{\text{d}}_H\mu=\beta$. Our proof above also shows that that the random limit measure $\nu$ fulfills $\CD\nu\geq 1-\alpha$ with conditional probability one. Moreover, an upper-bound technique due to Dawson and Hochberg (see \cite{DH,Zcharater}) can easily be applied in our setting to show that also the reverse inequality holds true. Thus, $$\PP(\CD\nu=1-\alpha|\nu\neq\text{null-measure})=1.$$
\end{remark}

\subsection{Proof of Theorem \ref{thm2}} This follows as in Theorem \ref{thm1} together from a modified Proposition \ref{prophelp}. The latter can be obtained by following the lines of the proofs of Lemmas 3.5 and 3.6 in \cite{BJST}. For these reasons the details are omitted.\hfill $\Box$

\subsection*{Acknowledgement} I am deeply indebted to Pierre Calka, who has drawn my attention to the problem considered in this note. Through his very inspiring talks at two workshops in Lille and Osnabr\"uck I had the privilege to learn from him about percolation in the hyperbolic plane.

\bibliographystyle{amsplain}

\begin{thebibliography}{20}

\bibitem{BenePetrino}
{R. Benedetti and C. Petrino}, {\em Lectures on Hyperbolic Geometry}, Springer, Berlin, 2008.

\bibitem{BJST}
{I. Benjamini, J. Jonasson, O. Schramm and J. Tykesson}, {\em Visibility to infinity in the hyperbolic plane, despite obstacles}, ALEA. Lat. Am. J. Probab. Math. Stat. \textbf{6} (2009), 323--342.

\bibitem{BenjaminiScramm}
{I. Benjamini and O. Schramm}, {\em Percolation in the hyperbolic plane}, J. Am. Math. Soc. \textbf{14} (2001), 487--507.

\bibitem{VisibilityPierre}
{P. Calka, J. Michel and  S. Porret-Blanc}, {\em Asymptotics of the visibility function in the Boolean model}, arXiv: 0905.4874 [math.PR] (2010).

\bibitem{CT}
{P. Calka and J. Tykesson}, {\em Asymptotics of visibility in the hyperbolic plane}, arXiv: 1012.5220 [math.PR] (2011).

\bibitem{Falconer}
{H. Falconer}, {\em Fractal Geometry. Mathematical Foundations and Applications}, 2nd Edition, Wiley, Chichester, 2003.

\bibitem{Federer}
{H. Federer}, {\em Geometric Measure Theory}, Springer, New York, 1969.

\bibitem{DH}
{D.A. Dawson and K.J. Hochberg}, {\em The carrying dimension of a stochastic measure diffusion}, Ann. Probab. \textbf{7} (1979), 693--703.

\bibitem{Kallenberg}
{O. Kallenberg}, {\em Random Measures}, 3rd edition Akademie Verlag, Berlin, 1983.

\bibitem{Kallenberg2}
{O. Kallenberg}, {\em Foundations of Modern Probability}, 2nd edition, Springer, New York, 2002.


\bibitem{Lalley}
{S. Lalley}, {\em Percolation clusters in hyperbolic tessellations}, Geom. Funct. Anal. \textbf{11} (2011), 971--1030.

\bibitem{MR}
{R. Meester and R. Roy}, {\em Continuum percolation}, Cambridge University Press, Cambridge, 1996.

\bibitem{Mol}
{I. Molchanov}, {\em Theory of Random Sets}, Springer, London, 2005.

\bibitem{Ramseyetal}
{A. Ramsay and R.D. Richtmyer}, {\em Introduction to Hyperbolic geometry}, Springer, New York, 2010.

\bibitem{Tyk}
{J. Tykesson}, {\em The number of unbounded components in the Poisson Boolean model of continuum percolation in hyperbolic space}, Electron. J. Probab. \textbf{12} (2007), 1379--1401.

\bibitem{ZCutout}
{U. Z\"ahle}, {\em Random fractals generated by random cutouts}, Math. Nachr. \textbf{116} (1984), 27--52.

\bibitem{Zcharater}
{U. Z\"ahle}, {\em The fractal character of localizable measure-valued processes, III. Fractal carrying sets of branching diffusions}, Math. Nachr. \textbf{138} (1988), 293--311.

\end{thebibliography}

\end{document}